\magnification=1200

\catcode`\À=\active \defÀ{\`A}    \catcode`\à=\active \defà{\`a} 
\catcode`\Â=\active \defÂ{\^A}    \catcode`\â=\active \defâ{\^a} 
\catcode`\Æ=\active \defÆ{\AE}    \catcode`\æ=\active \defæ{\ae}
\catcode`\Ç=\active \defÇ{\c C}   \catcode`\ç=\active \defç{\c c}
\catcode`\È=\active \defÈ{\`E}    \catcode`\è=\active \defè{\`e} 
\catcode`\É=\active \defÉ{\'E}    \catcode`\é=\active \defé{\'e} 
\catcode`\Ê=\active \defÊ{\^E}    \catcode`\ê=\active \defê{\^e} 
\catcode`\Ë=\active \defË{\"E}    \catcode`\ë=\active \defë{\"e} 
\catcode`\Î=\active \defÎ{\^I}    \catcode`\î=\active \defî{\^\i}
\catcode`\Ï=\active \defÏ{\"I}    \catcode`\ï=\active \defï{\"\i}
\catcode`\Ô=\active \defÔ{\^O}    \catcode`\ô=\active \defô{\^o} 
\catcode`\Ù=\active \defÙ{\`U}    \catcode`\ù=\active \defù{\`u} 
\catcode`\Û=\active \defÛ{\^U}    \catcode`\û=\active \defû{\^u} 
\catcode`\Ü=\active \defÜ{\"U}    \catcode`\ü=\active \defü{\"u} 

\catcode`\ =\active \def { }

\hsize=11.25cm    
\vsize=18cm       
\parindent=12pt   \parskip=5pt     

\hoffset=.5cm   
\voffset=.8cm   

\pretolerance=500 \tolerance=1000  \brokenpenalty=5000

\catcode`\@=11

\font\eightrm=cmr8         \font\eighti=cmmi8
\font\eightsy=cmsy8        \font\eightbf=cmbx8
\font\eighttt=cmtt8        \font\eightit=cmti8
\font\eightsl=cmsl8        \font\sixrm=cmr6
\font\sixi=cmmi6           \font\sixsy=cmsy6
\font\sixbf=cmbx6

\font\tengoth=eufm10 
\font\eightgoth=eufm8  
\font\sevengoth=eufm7      
\font\sixgoth=eufm6        \font\fivegoth=eufm5

\skewchar\eighti='177 \skewchar\sixi='177
\skewchar\eightsy='60 \skewchar\sixsy='60

\newfam\gothfam           \newfam\bboardfam

\def\tenpoint{
  \textfont0=\tenrm \scriptfont0=\sevenrm \scriptscriptfont0=\fiverm
  \def\rm{\fam\z@\tenrm}
  \textfont1=\teni  \scriptfont1=\seveni  \scriptscriptfont1=\fivei
  \def\oldstyle{\fam\@ne\teni}\let\old=\oldstyle
  \textfont2=\tensy \scriptfont2=\sevensy \scriptscriptfont2=\fivesy
  \textfont\gothfam=\tengoth \scriptfont\gothfam=\sevengoth
  \scriptscriptfont\gothfam=\fivegoth
  \def\goth{\fam\gothfam\tengoth}
  
  \textfont\itfam=\tenit
  \def\it{\fam\itfam\tenit}
  \textfont\slfam=\tensl
  \def\sl{\fam\slfam\tensl}
  \textfont\bffam=\tenbf \scriptfont\bffam=\sevenbf
  \scriptscriptfont\bffam=\fivebf
  \def\bf{\fam\bffam\tenbf}
  \textfont\ttfam=\tentt
  \def\tt{\fam\ttfam\tentt}
  \abovedisplayskip=12pt plus 3pt minus 9pt
  \belowdisplayskip=\abovedisplayskip
  \abovedisplayshortskip=0pt plus 3pt
  \belowdisplayshortskip=4pt plus 3pt 
  \smallskipamount=3pt plus 1pt minus 1pt
  \medskipamount=6pt plus 2pt minus 2pt
  \bigskipamount=12pt plus 4pt minus 4pt
  \normalbaselineskip=12pt
  \setbox\strutbox=\hbox{\vrule height8.5pt depth3.5pt width0pt}
  \let\bigf@nt=\tenrm       \let\smallf@nt=\sevenrm
  \normalbaselines\rm}

\def\eightpoint{
  \textfont0=\eightrm \scriptfont0=\sixrm \scriptscriptfont0=\fiverm
  \def\rm{\fam\z@\eightrm}
  \textfont1=\eighti  \scriptfont1=\sixi  \scriptscriptfont1=\fivei
  \def\oldstyle{\fam\@ne\eighti}\let\old=\oldstyle
  \textfont2=\eightsy \scriptfont2=\sixsy \scriptscriptfont2=\fivesy
  \textfont\gothfam=\eightgoth \scriptfont\gothfam=\sixgoth
  \scriptscriptfont\gothfam=\fivegoth
  \def\goth{\fam\gothfam\eightgoth}
  
  \textfont\itfam=\eightit
  \def\it{\fam\itfam\eightit}
  \textfont\slfam=\eightsl
  \def\sl{\fam\slfam\eightsl}
  \textfont\bffam=\eightbf \scriptfont\bffam=\sixbf
  \scriptscriptfont\bffam=\fivebf
  \def\bf{\fam\bffam\eightbf}
  \textfont\ttfam=\eighttt
  \def\tt{\fam\ttfam\eighttt}
  \abovedisplayskip=9pt plus 3pt minus 9pt
  \belowdisplayskip=\abovedisplayskip
  \abovedisplayshortskip=0pt plus 3pt
  \belowdisplayshortskip=3pt plus 3pt 
  \smallskipamount=2pt plus 1pt minus 1pt
  \medskipamount=4pt plus 2pt minus 1pt
  \bigskipamount=9pt plus 3pt minus 3pt
  \normalbaselineskip=9pt
  \setbox\strutbox=\hbox{\vrule height7pt depth2pt width0pt}
  \let\bigf@nt=\eightrm     \let\smallf@nt=\sixrm
  \normalbaselines\rm}

\tenpoint

\def\pc#1{\bigf@nt#1\smallf@nt}         \def\pd#1 {{\pc#1} }

\catcode`\;=\active
\def;{\relax\ifhmode\ifdim\lastskip>\z@\unskip\fi
\kern\fontdimen2  -1.2 \fontdimen3 \string;}

\catcode`\:=\active
\def:{\relax\ifhmode\ifdim\lastskip>\z@\unskip\fi\penalty\@M\ \fi\string:}

\catcode`\!=\active
\def!{\relax\ifhmode\ifdim\lastskip>\z@
\unskip\fi\kern\fontdimen2  -1.1 \fontdimen3 \string!}

\catcode`\?=\active
\def?{\relax\ifhmode\ifdim\lastskip>\z@
\unskip\fi\kern\fontdimen2  -1.1 \fontdimen3 \string?}

\catcode`\«=\active 
\def«{\raise.4ex\hbox{%
 $\scriptscriptstyle\langle\!\langle$}}

\catcode`\»=\active 
\def»{\raise.4ex\hbox{%
 $\scriptscriptstyle\rangle\!\rangle$}}

\frenchspacing

\def\raggedbottom{\topskip 10pt plus 36pt\r@ggedbottomtrue}

\def\pointir{\unskip . --- \ignorespaces}

\def\Medbreak{\vskip-\lastskip\medbreak}

\long\def\th#1 #2\enonce#3\endth{
   \Medbreak\noindent
   {\pc#1} {#2\unskip}\pointir{\it #3}\smallskip}

\def\decale#1{\smallbreak\hskip 28pt\llap{#1}\kern 5pt}
\def\decaledecale#1{\smallbreak\hskip 34pt\llap{#1}\kern 5pt}
\def\puce{\smallbreak\hskip 6pt{$\scriptstyle\bullet$}\kern 5pt}

\def\eqalign#1{\null\,\vcenter{\openup\jot\m@th\ialign{
\strut\hfil$\displaystyle{##}$&$\displaystyle{{}##}$\hfil
&&\quad\strut\hfil$\displaystyle{##}$&$\displaystyle{{}##}$\hfil
\crcr#1\crcr}}\,}

\catcode`\@=12

\showboxbreadth=-1  \showboxdepth=-1

\newcount\numerodesection \numerodesection=1
\def\section#1{\bigbreak
 {\bf\number\numerodesection.\ \ #1}\nobreak\medskip
 \advance\numerodesection by1}

\mathcode`A="7041 \mathcode`B="7042 \mathcode`C="7043 \mathcode`D="7044
\mathcode`E="7045 \mathcode`F="7046 \mathcode`G="7047 \mathcode`H="7048
\mathcode`I="7049 \mathcode`J="704A \mathcode`K="704B \mathcode`L="704C
\mathcode`M="704D \mathcode`N="704E \mathcode`O="704F \mathcode`P="7050
\mathcode`Q="7051 \mathcode`R="7052 \mathcode`S="7053 \mathcode`T="7054
\mathcode`U="7055 \mathcode`V="7056 \mathcode`W="7057 \mathcode`X="7058
\mathcode`Y="7059 \mathcode`Z="705A


\catcode`\À=\active \defÀ{\`A}    \catcode`\à=\active \defà{\`a} 
\catcode`\Â=\active \defÂ{\^A}    \catcode`\â=\active \defâ{\^a} 
\catcode`\Æ=\active \defÆ{\AE}    \catcode`\æ=\active \defæ{\ae}
\catcode`\Ç=\active \defÇ{\c C}   \catcode`\ç=\active \defç{\c c}
\catcode`\È=\active \defÈ{\`E}    \catcode`\è=\active \defè{\`e} 
\catcode`\É=\active \defÉ{\'E}    \catcode`\é=\active \defé{\'e} 
\catcode`\Ê=\active \defÊ{\^E}    \catcode`\ê=\active \defê{\^e} 
\catcode`\Ë=\active \defË{\"E}    \catcode`\ë=\active \defë{\"e} 
\catcode`\Î=\active \defÎ{\^I}    \catcode`\î=\active \defî{\^\i}
\catcode`\Ï=\active \defÏ{\"I}    \catcode`\ï=\active \defï{\"\i}
\catcode`\Ô=\active \defÔ{\^O}    \catcode`\ô=\active \defô{\^o} 
\catcode`\Ù=\active \defÙ{\`U}    \catcode`\ù=\active \defù{\`u} 
\catcode`\Û=\active \defÛ{\^U}    \catcode`\û=\active \defû{\^u} 
\catcode`\Ü=\active \defÜ{\"U}    \catcode`\ü=\active \defü{\"u}

\def\ogoth{\goth o}

\centerline{\bf Good reduction, bad reduction}
\bigskip
\centerline{\it Notes for a lecture at the conference on}
\smallskip
\centerline{\it Commutative Algbera and Algebraic Geometry, Madras,
August 1--6, 2005.}
\bigskip
We give some general properties of good and bad reduction, along with
some recent examples (worked out with Dipendra Prasad) of varieties
having bad reduction not accounted for by cohomology.  We include some
consequences of our remarks for varieties over number fields having
good reduction everywhere.
\bigskip
\section{Good reduction.}

Let $p$ be a prime number and let $K$ be a finite extension of
${\bf Q}_p$.  
Denote by $\ogoth$ the ring of integers of $K$
and let $k$ be the residue field of
$\ogoth$, the quotient by the unique maximal ideal.

\th DEFINITION 
\enonce
A smooth proper $K$-variety $X$ is said to have good reduction if it is
the generic fibre of a smooth proper $\ogoth$-scheme.
\endth

Such an $\ogoth$-scheme is called a smooth model of $X$.  A variety is
said to have {\it bad reduction\/} if it does not have good reduction.

As an example, a finite extension $L$ of $K$ has good reduction if and
only if it is unramified.


There are exactly two conics (i.e.~curves of genus~$0$) over $K$.  The
one which has a point, namely ${\bf P}_1$, has good reduction; the
other has bad reduction.

More generally, a twisted form of ${\bf P}_n$ has bad reduction unless
it has a $K$-point, and is thus isomorphic to ${\bf P}_n$.

A twisted form of an abelian variety has bad reduction if it does not
have a $K$-point (it may have bad reduction even when it has one).

An elliptic curve has good reduction if and only if the its conductor
is (the unit ideal) $\ogoth$.


There are some varieties which have (infinitely) many smooth models :
A smooth model of ${\bf P}_1\times {\bf P}_1$ can be found whose
special fibre is any of the the Hirzebruch surfaces $F_0={\bf
P}_1\times {\bf P}_1$, $F_2$, $F_4$, $\ldots$ I don't know of a
variety which has at least two but only finitely many smooth models.

Let $l$ be a prime number different from $p$.  Let $\bar K$ be an
algebraic closure of $K$.  For a $K$-variety $X$, denote by $\bar
X=X\times_K\bar K$ the change of base field from $K$ to $\bar K$.

\th THEOREM
\enonce
Let $X$ be a (smooth, proper) $K$-variety.  If $X$ has good reduction,
then the action of ${\rm Gal}(\bar K|K)$ on $H^i(\bar X, {\bf Q}_l)$
is unramified.
\endth


This means that the action factors through the quotient ${\rm
Gal}(\tilde K|K)$, where $\tilde K$ is the maximal unramified
extension of $K$ in $\bar K$ : the inertia subgroup ${\rm Gal}(\bar
K|\tilde K)$ acts trivially.

$p$-adic chomology is almost never unramified.  It took Fontaine to
formulate the correct analogue, and Faltings to prove it in general.

\th THEOREM 
\enonce
Let $X$ be a (smooth, proper) $K$-variety.  If $X$ has good reduction,
then the action of ${\rm Gal}(\bar K|K)$ on $H^i(\bar X, {\bf Q}_p)$
is crystalline.
\endth

This is also a condition on the restriction of the representation to
the inertial subgroup.  More precisely, Fontaine has constructed a
${\bf Q}_p$-algebra $B_{\rm cris}$ with an action of ${\rm Gal}(\bar
K|K)$ (and some other structures) ; a finite-dimensional ${\bf
Q}_p$-representation $V$ of ${\rm Gal}(\bar K|K)$ is called
crystalline if the dimension of $(V\otimes_{{\bf Q}_p}B_{\rm
cris})^{{\rm Gal}(\bar K|K)}$ (as a vector space over $B_{\rm
cris}^{{\rm Gal}(\bar K|K)}$, the maximal unramified extension of
${\bf Q}_p$ in $K$) is the same as that of $V$.


\th QUESTION
\enonce
Is there a $K$-variety $X$ for which the action of ${\rm Gal}(\bar
K|K)$ on $H^i(\bar X, {\bf Q}_l)$ is unramified, and the action on
$H^i(\bar X, {\bf Q}_p)$ crystalline, and yet which has bad reduction.
\endth

%
%
%

We shall provide some examples of such varieties, which have bad
reduction but whose ``motive'' has good reduction.

\bigbreak
\centerline{\bf Twisted abelian varieties}

\th DEFINITION
\enonce
An abelian $K$-variety $A$ is said to have abelian reduction if
there is an abelian $\ogoth$-scheme whose generic fibre is $A$.
\endth

The abelian scheme in question is then unique.  Clearly, if $A$ has
abelian reduction, then it has good reduction.  Conversely,

\th THEOREM
\enonce
If the variety $A$ has good reduction, then the abelian
variety $A$ has abelian reduction.
\endth

A direct proof can be found in the book on Néron Models by Bosch,
L{\"u}tkebohmert Raynaud.  

%
%
%
%
%

Here is the celebrated $l$-adic criterion (Néron-Ogg-Shafarevich) for
good reduction of abelian varieties, proved by Serre and Tate.

\th THEOREM 
\enonce
Let $A$ be an abelian $K$-variety.  If the representation $H^1(\bar X,
{\bf Q}_l)$ is unramified, then $A$ has good reduction.
\endth

So if the variety $A$ has good reduction (disregarding the group law),
the $l$-adic cohomology is unramified and hence $A$ has abelian
reduction.

The $p$-adic analogue is more recent, and due to Mokrane and Coleman-Iovita.

\th THEOREM 
\enonce
Let $A$ be an abelian $K$-variety.  If the representation $H^1(\bar X,
{\bf Q}_p)$ is crystalline, then $A$ has good reduction.
\endth

Our observation for abelian varieties amounts to the next two
results. 

\th PROPOSITION
\enonce
Let\/ $T$ be a torsor under an abelian variety $A$.  Then the
representation of\/ ${\rm Gal}(\bar K|K)$ on the (\/$l$-adic or
$p$-adic) cohomology of\/ $T$ is the same as the representation on the
cohomology of $A$.
\endth

In short, the motive of $T$ is the same as that of $A$.


\th PROPOSITION
\enonce
Let\/ $T$ be a $K$-variety which is potentially isomorphic to an abelian
variety.  If\/ $T(K)$ is empty, then\/ $T$ has bad reduction.
\endth

Thus taking an $A$ which has good reduction and a torsor $T$ which is
not $K$-isomorphic to $A$ (there are many such $T$, by Tate local
duality), the variety $T$ has bad reduction but its cohomology is
unramified (resp. crystalline).

\bigbreak
\centerline{\bf Twisted projective spaces}

%

These are the varieties which become isomorphic to ${\bf P}_n$ over a
suitable finite (separable) extension of $K$.  They are classified by
the group $H^2(K,{\bar K}^\times)$, as are similarity classes of
simple central $K$-algebras.

\th PROPOSITION
\enonce
Let\/ $X$ be a twisted form of\/ ${\bf P}_n$ and\/ $A$ the corresponding
simple central\/ $K$-algebra.  Then the following statements are
equivalent :

1) $X$ is $K$-isomorphic to ${\bf P}_n$.

2) $X$ has good reduction.

3) $X(K)$ is not empty.


4) $A$ is similar to the matrix algebra.

5) $A$ is the generic fibre of an azumaya $\ogoth$-algebra.

6) the class of\/ $X$ in $H^2(K,{\bar K}^\times)$ is trivial.

\endth

As in the earlier case of abelian varieties, the cohomology of a
twisted form is the twist by (the image of) the same $1$-cocycle.
This gives :

\th PROPOSITION
\enonce
The $l$-adic (resp.~$p$-adic) cohomology of a twisted form $X$ of
${\bf P}_n$ is unramified (resp.~crystalline).
\endth

Taking $X$ to be different from ${\bf P}_n$, we get a variety which
has bad reduction but whose cohomology is unramified
(resp.~crystalline).  Recall that $H^2(K,{\bar K}^\times)={\bf Q}/{\bf
Z}$, so there are many such varieties.

\bigbreak
\centerline{\bf Rational surfaces}

These are the surfaces which are potentially (i.e.~over a suitable
finite (separable) extension of~$K$) birational to ${\bf P}_2$.  

Let us confine ourselves to Ch{\^a}telet surfaces and to $p\neq2$.

Let $d\in K^\times$ not be a square and $e_1, e_2$ be two distinct
elements of $K^\times$.  The ruled
surface $X$ given by
$$
y^2-dz^2=xx'(x-e_1x')(x-e_2x')t^2
$$ 
is fibered in conics over ${\bf P}_1$ (coordinates $x:x'$), the fibre
at each point being a conic in ${\bf P}_2$ (coordinates $y:z:t$).  The
surface is birational over $K(\sqrt d)$ to ${\bf P}_2$.

Without changing the surface $X$, we may assume that $e_1, e_2$ have
the same valuation $r$.

%
%
%
%
%
%

\th PROPOSITION
\enonce
The surface $X$ has bad reduction if\/ $d$ is not a unit, or if
the valuation of\/ $e_1-e_2$ is $>r$. 
\endth
The reason is that in these cases the Chow group of $0$-cycles of
degree~$0$ turns out to be $\neq0$, whereas a theorem of
Colliot-Th{\'e}l{\`e}ne says that it should vanish for rational
surfaces having good reduction.

\th PROPOSITION
\enonce
If\/ $d$ is a unit (and $p\neq2$), the $l$-adic (resp.~$p$-adic)
cohomology of\/ $X$ is unramified (resp.~crystalline).
\endth
This follows from the fact that the cohomology of a rational surface
can be computed from its Picard group over $\bar K$ (a finitely
generated free ${\bf Z}$-module with continuous ${\rm Gal}(\bar
K|K)$-action).

Taking $d$ to be a unit and the valuation of $e_1-e_2$ to be $>r$, we
get examples of rational surfaces which have bad reduction but whose
``motive'' has good reduction.

\bigbreak
\centerline{\bf Curves of higher genus}

The most natural way of finding examples among curves of genus $\ge2$
would be to combine N{\'e}ron-Ogg-Shafarevich with the anabelian
$l$-adic criterion of Oda :

\th THEOREM
\enonce
Let\/ $C$ be a curve of genus $\ge2$.  If the outer action of\/ ${\rm
Gal}(\bar K|K)$ on the maximal pro-$l$ quotient of the ({\'e}tale)
fundamental group of\/ $C_{\bar K}=C\times_K\bar K$ is unramified,
then $C$ has good reduction.
\endth

%

\bigbreak
\centerline{\bf Over number fields}

The classical result of Minkowski about unramified extensions of ${\bf
Q}$ can be reformulated as :
\th THEOREM
\enonce
The only point over ${\bf Q}$ which has good reduction
everywhere is\/ ${\rm Spec}({\bf Q})$.
\endth

In dimension~$1$, there is something similar :
\th THEOREM
\enonce
The only (smooth, proper, absolutely connected) curve over ${\bf Q}$
which has good reduction everywhere is\/ ${\bf P}_1$.
\endth

Our observations suffice to prove this in genus~$0$ and, combined with
the theorem of Tate saying that there is no elliptic curve over ${\bf
Q}$ having good reduction everywhere, in genus~$1$.  In higher genera,
one uses the fact that if a curve has good reduction, then so does its
jacobian --- an easy consequence of N{\'e}ron-Ogg-Shafarevich, but
there is also a direct proof --- and the famous theorem of Fontaine,
generalising Tate's result to higher dimensions :

\th THEOREM
\enonce
The only abelian variety over ${\bf Q}$ which has good
reduction everywhere is the point.
\endth

These two theorems, about smooth points (resp. curves) over ${\bf Z}$,
allow N. Fakhruddin to show that the degree-4 (resp.~degree-3)
embedding ${\bf P}_1\rightarrow{\bf P}_4$ (resp.~${\bf
P}_2\rightarrow{\bf P}_9$) does not admit a smooth hypersurface
section (cf.~Poonen's recent Annals paper).


Our results about twisted forms of projective spaces imply :

\th THEOREM
\enonce
The only twisted form of\/ ${\bf P}_n$ over\/ ${\bf Q}$ which has good
reduction everywhere is\/ ${\bf P}_n$.
\endth

Over a given number field, the possible twisted forms of ${\bf P}_n$
having good reduction everywhere can easily be listed, basically
because there is a local-to-global principle for them
(Brauer-Hasse-N{\oe}ther).  For example, over a real quadratic field,
there is a unique curve of genus~$0$, apart from ${\bf P}_1$, which
has good reduction everywhere.
\vfill
\line{\hfill Chandan Singh Dalawat}
\line{\hfill Harish-Chandra Research Institute}
\line{\hfill Chhatnag Road, Jhunsi}
\line{\hfill Allahabad 211 019}

\bye